\documentclass[12pt,reqno,a4wide]{amsart}

\usepackage{tikz} 
\usepackage{calrsfs}
\usepackage{mathrsfs}
\usetikzlibrary{decorations.pathmorphing}
\usetikzlibrary{arrows.meta}
\usepgflibrary {shadings}
\usepackage{array}
\usepackage{xcolor}
\usepackage{comment}
\usepackage{hyperref}

\usetikzlibrary{decorations.pathreplacing}
\usetikzlibrary{fit}

\usepackage{pgfplots}

\allowdisplaybreaks

    \def\ss{\sigma}

\catcode`,\active

\catcode`\,12

\begin{document}
\bibliographystyle{plain}

%
%

	\title
	{Arndt-Carlitz compositions}

	\author[H. Prodinger ]{Helmut Prodinger }
	\address{Department of Mathematics, University of Stellenbosch 7602, Stellenbosch, South Africa
	and
NITheCS (National Institute for
Theoretical and Computational Sciences), South Africa.}
	\email{hproding@sun.ac.za}


	\begin{abstract}
Carlitz-compositions follow the restrictions of neighbouring parts $\ss_{i-1}\neq\ss_{i}$. The recently introduced Arndt-compositions have to satisfy
$\ss_{2i-1}>\ss_{2i}$. The two concepts are combined to new and exciting objects that we call Arndt-Carlitz compositions.

	\end{abstract}
	

\maketitle

\section{Introduction}
A composition of a (positive) integer $n$ is a representation $n=\ss_1+\dots+\ss_k$ with positive integers $\ss_i$; $k$ is called the number of parts.
Compositions with (local) restrictions have been popular for many years. For instance, some 25 years ago, Knopfmacher and Prodinger investigated so-called Carlitz compositions \cite{KP}, where neighbouring parts have to be different.

Recently, another (simpler) concept became popular: Arndt-compositions. Here, $\ss_{2i-1}>\ss_{2i}$ for all possible indices $i$ \cite{HT1,HT2,Ramirez}. The aim of the present  note is to combine both concepts and consider what we call Arndt-Carlitz compositions. We only concentrate on basic enumeration of these objects, including asymptotics and leave various related quantities to other interested researchers.

\section{Arndt-Carlitz compositions with an even number of parts}

The restrictions we have to deal with are
\begin{equation*}
\ss_1>\ss_2\neq\ss_3>\ss_4\neq \ss_5\dots\ss_{2k-1}>\ss_{2k}.
\end{equation*}

We compute a bivariate generating function, as in \cite{KP} and use the adding-a-new slice technique, see \cite{FP}.
\begin{equation*}
F(z,u)=		\sum_{k\ge1}a_k(z,u);
\end{equation*}
the variable $z$ is used to record the composition, and $u$ to record the value of $\ss_{2k}$. The variable $k$ is used to count Arndt-Carlitz compositions with
$2k$ parts.  

Adding a new slice means here to add the new parts, satisfying the conditions. The exponent of $u$ is essential here:
\begin{align*}
u^j&\longrightarrow\sum_{i\neq j}\sum_{k<i}z^{i+k}u^k=\sum_{i}\sum_{k<i}z^{i+k}u^k-\sum_{k<j}z^{j+k}u^k\\&=
\frac{z^3u}{(1-z)(1-z^2u)}-\frac{z^j(zu-(zu)^{j})}{1-zu}.
\end{align*}
Then we get $a_{k+1}(z,u)$ from $a_{k}(z,u)$:
\begin{align*}
	a_{k+1}(z,u)&=
	\frac{z^3u}{(1-z)(1-z^2u)}a_k(z,1)-\frac{zu}{1-zu}a_k(z,z)+\frac{1}{1-zu}a_k(z,z^2u),\\
 a_1(z,u)&=\frac{z^3u}{(1-z)(1-z^2u)}.
\end{align*}
Summing, we get the essential functional equation
\begin{align*}
F&(z,u)=		\sum_{k\ge1}a_k(z,u)=a_1(z,u)+\sum_{k\ge1}a_{k+1}(z,u)\\*
&=\sum_{k\ge1}\bigg[	\frac{z^3u}{(1-z)(1-z^2u)}a_k(z,1)-\frac{zu}{1-zu}a_k(z,z)+\frac{1}{1-zu}a_k(z,z^2u)\biggr]\\&\quad+\frac{z^3u}{(1-z)(1-z^2u)}\\
&=\frac{z^3u}{(1-z)(1-z^2u)}F(z,1)-\frac{zu}{1-zu}F(z,z)+\frac{1}{1-zu}F(z,z^2u)+\frac{z^3u}{(1-z)(1-z^2u)}.
\end{align*}
This will be iterated:
\begin{align*}
	F(z,u)
&=\frac{z^3u}{(1-z)(1-z^2u)}F(z,1)-\frac{zu}{1-zu}F(z,z)+\frac{z^3u}{(1-z)(1-z^2u)}\\
&+\frac{1}{1-zu}\biggl[\frac{z^5u}{(1-z)(1-z^4u)}F(z,1)-\frac{z^3u}{1-z^3u}F(z,z)+\frac{z^5u}{(1-z)(1-z^4u)}\\
&+\frac{1}{1-z^3u}\biggl[\frac{z^7u}{(1-z)(1-z^6u)}F(z,1)-\frac{z^5u}{1-z^5u}F(z,z)+\frac{z^7u}{(1-z)(1-z^6u)}\\&
+\frac{1}{1-z^5u}\biggl[\dots
\end{align*}
Collecting,
\begin{align*}
	F(z,u)
	&=\frac{z^3u}{(1-z)(1-z^2u)}F(z,1)+\frac{z^5u}{(1-z)(1-zu)(1-z^4u)}F(z,1)\\&+\frac{z^7u}{(1-z)(1-zu)(1-z^3u)(1-z^6u)}F(z,1)+\dots\\
	&-\frac{zu}{1-zu}F(z,z)-\frac{z^3u}{(1-zu)(1-z^3u)}F(z,z)\\&-\frac{z^5u}{(1-zu)(1-z^3u)(1-z^5u)}F(z,z)-\dots\\
	&+\frac{z^3u}{(1-z)(1-z^2u)}+\frac{1}{1-zu}\frac{z^5u}{(1-z)(1-z^4u)}+\frac{1}{(1-zu)(1-z^3u)}\frac{z^7u}{(1-z)(1-z^6u)}+\dots\\	
	 \end{align*}
 or
 \begin{align*}
 	F(z,u)
 	&=\frac{1}{1-z}\biggl[\frac{z^3u}{(1-z^2u)}+\frac{z^5u}{(1-zu)(1-z^4u)}\\&\qquad\qquad\qquad+\frac{z^7u}{(1-zu)(1-z^3u)(1-z^6u)}+\dots\biggr](F(z,1)+1)\\
 	&-\biggl[\frac{zu}{1-zu}+\frac{z^3u}{(1-zu)(1-z^3u)}+\frac{z^5u}{(1-zu)(1-z^3u)(1-z^5u)}+\dots\biggr]F(z,z),
 \end{align*}
or
\begin{equation*}
F(z,u)=\alpha(z,u)(F(z,1)+1)+\beta(z,u)F(z,z).
\end{equation*}
We are (in this section) only interested in $F(z,1)$, and we find the system 
\begin{align*}
F(z,1)&=\alpha(z,1)(F(z,1)+1)+\beta(z,1)F(z,z),\\
F(z,z)&=\alpha(z,z)(F(z,1)+1)+\beta(z,z)F(z,z),
\end{align*}
with
\begin{gather*}
\alpha(z,u)=\frac1{1-z}\sum_{k\ge1}\frac{z^{2k+1}u}{1-z^{2k}u}\dfrac1{\prod_{\ell=1}^{k-1}(1-z^{2\ell-1}u)},\\
\beta(z,u)=-\sum_{k\ge1}\dfrac{z^{2k-1}u}{\prod_{\ell=1}^k(1-z^{2\ell-1}u)}.
\end{gather*}
The special values are then
\begin{equation*}
\alpha(z,1)=\frac{z}{1-z}\sum_{k\ge1}\frac{z^{2k}}{1-z^{2k}}\dfrac1{\prod_{\ell=1}^{k-1}(1-z^{2\ell-1})},\quad
\alpha(z,z)=\frac z{1-z}\sum_{k\ge1}\frac{z^{2k+1}}{1-z^{2k+1}}\dfrac1{\prod_{\ell=1}^{k-1}(1-z^{2\ell})}
\end{equation*}
and
\begin{equation*}
	\beta(z,1)=-\sum_{k\ge1}\dfrac{z^{2k-1}}{\prod_{\ell=1}^k(1-z^{2\ell-1})},\quad
	\beta(z,z)=-\sum_{k\ge1}\dfrac{z^{2k}}{\prod_{\ell=1}^k(1-z^{2\ell})}.
\end{equation*}
Of course, once $F(z,1)$ and $F(z,z)$ are known, $F(z,u)$ is also known.

Solving the system of two equations leads to the final answer
\begin{equation*}
F(z,1)=\frac{\alpha(z,1)+\alpha(z,z)\beta(z,1)-\alpha(z,1)\beta(z,z)}{1-\alpha(z,1)-\beta(z,z)+\beta(z,z)\alpha(z,1)-\alpha(z,z)\beta(z,1)}
\end{equation*}
and the series expansion
\begin{equation*}
F(z,1)=z^3+z^4+2z^5+3z^6+5z^7+7z^8+12z^9+20z^{10}+30z^{11}+\dots
\end{equation*}
Here are the 5 Arndt-Carlitz compositions of 7 into an even number of parts:
$6+1,5+2,4+3,3+1+2+1,2+1+3+1$ and the 7 compositions of 8: $7+1,6+2,5+3,3+1+3+1,2+1+4+1,2+1+3+2,4+1+2+1$.

\section{Arndt-Carlitz compositions with an odd number of parts}

First, we need the other solution of the system,
\begin{equation*}
F(z,z)=\frac{\alpha(z,z)}{1-\alpha(z,1)-\beta(z,z)+\beta(z,z)\alpha(z,1)-\alpha(z,z)\beta(z,1)}.
\end{equation*}
Now we get such a composition (with an odd number of parts) by the compositions from the previous section (even number of parts) and attaching an extra number, which must be different from previous last number. This leads to
\begin{equation*}
\frac{z}{1-z}+F(z,1)\frac{z}{1-z}-F(z,z).
\end{equation*}
The first term is responsible for compositions into just one part, which we would not get by the attaching procedure.
The series expansion is 
\begin{equation*}
\frac{z}{1-z}+F(z,1)\frac{z}{1-z}-F(z,z)=z+z^2+z^3+z^4+2z^5+4z^6+5z^7+9z^8+15z^9+22z^{10}+36z^{11}+\dots
\end{equation*}
Here are the 9 Arndt-Carlitz compositions of 8 into an odd number of parts: $8,2+1+5,3+1+4,3+2+3,4+1+3,5+1+2,5+2+1,4+3+1,2+1+2+1+2$.

\section{Asymptotics}

The strategy to do this is exactly as in \cite{FP} and \cite{KP}; there is a simple pole of the respective generating functions, which is closest to the origin.
It can be computed numerically. The justification that the procedure works well is exactly as in the older papers and shall not be repeated here.
The numerator of $F(z,1)$ has its closest zero to the origin at
$\rho=.62790101012637517122$. Then
\begin{equation*}
F(z,1)\sim 0.18236796484521070938\frac1{1-z/\rho};
\end{equation*}
\begin{equation*}
\frac{z}{1-z}+F(z,1)\frac{z}{1-z}-F(z,z)\sim0.217010508476828474\frac1{1-z/\rho}
\end{equation*}
Adding the constants, $0.18236796484521070938+0.217010508476828474=0.399378473322039$, and
$\rho^{-1}=1.592607726174439$, and so we have an asymptotic formula for the Arndt-Carlitz compositions of $n$ (no parity restrictions on the number of parts):
\begin{equation*}
0.399378473322039\cdot 1.592607726174439^n.
\end{equation*}
For unrestricted compositions the formula is $\textsf{const}\cdot2^n$ and for Carlitz-compositions $\textsf{const}\cdot1.750243^n$. 
That we get a smaller number $1.5926$ in the instance of Arndt-Carlitz compositions is not surprising since the restrictions are a bit more strict.
\bibliographystyle{plain}

\end{document}